\documentclass{amsart}

\newtheorem{theorem}{Theorem}[section]
\newtheorem{lemma}[theorem]{Lemma}

\theoremstyle{definition}

\newtheorem{example}[theorem]{Example}

\newtheorem{prop}[theorem]{Proposition}
\theoremstyle{remark}
\newtheorem{remark}[theorem]{Remark}
\usepackage{mathrsfs} \numberwithin{equation}{section}
\newcommand{\R}{\mathbb{R}}
\newcommand{\C}{\mathbb{C}}
\newcommand{\Z}{\mathbb{Z}}
\newcommand{\N}{\mathbb{N}}

\newcommand{\PP}{\mathbb{P}}
\begin{document}

\title{
Application of multivariate splines to discrete mathematics}

%    Information for first author
\author{Zhiqiang Xu}
%    Address of record for the research reported here
\address{Institute of Computational Mathematics, Academy of Mathematics and
System Sciences, Chinese Academy of Sciences, Beijing, 100080,
China}
%    Current address
%\curraddr{Department of Mathematics and Statistics, Case Western
%Reserve University, Cleveland, Ohio 43403}
\email{xuzq@lsec.cc.ac.cn, Fax: +86-10-62542285.}
%    \thanks will become a 1st page footnote.

%    Information for second author

%    General info

%\date{January 1, 2001 and, in revised form, June 22, 2001.}

%\dedicatory{This paper is dedicated to our advisors.}


\begin{abstract}
Using  methods developed in multivariate splines, we present an
explicit formula for discrete truncated powers, which are defined
as the number of non-negative integer solutions of linear
Diophantine equations. We further use the formula to study some
classical problems in discrete mathematics as follows. First, we
extend the partition function of integers in number theory.
Second, we exploit the relation between the relative volume of
convex polytopes and multivariate truncated powers and give a
simple proof for the volume formula for the Pitman-Stanley
polytope. Third, an explicit formula for the Ehrhart
quasi-polynomial is presented.
\end{abstract}
\baselineskip 20pt

 \maketitle
\section{Introduction}

\setcounter{section}{1}

 Let $M$ be an $s\times n$ integer matrix
with columns $m_1,\ldots,m_n\in {\Z}^s$ such that their convex
hull does not contain the origin. For a given $\alpha \in {\Z}^s$,
consider the following system of linear Diophantine equations
$$
M\beta\,\,=\,\,\alpha,\,\,\, \beta\in {\Z}^n.
$$
The number of non-negative integer solutions $\beta$ for this system
is denoted by $t(\alpha |M)$ and the resulting function $t(\cdot|M)$
on $\Z^s$ is called a {\it discrete truncated power}. Discrete
truncated powers,  also called {\it vector partition functions},
have many applications in various mathematical areas including
Algebraic Geometry, Representation Theory, Number Theory, Statistics
and Randomized Algorithms. In general, one studies $t(\cdot |M)$ by
generating functions and by means of algebraic geometry etc.
\cite{BeckFrobenius,beckbook,brion,cappell,vergne1}. For example,
when $M$ is unimodular, that is when every nonsingular square
submatrix of order $s$ has determinant $\pm 1$, two algebraic
algorithms for generating the explicit formula for $t({\cdot }|M)$
are presented in \cite{deloera}. When $s=1,$ an explicit formula for
$t({\cdot }|M)$, which counts the integer solutions for the linear
Diophantine equation, is presented in \cite{BeckFrobenius}.
Especially,  Popoviciu  \cite{popoviciu} gave a beautiful formula
for $t(\cdot |M)$, when $M=(a, b)$ where $a$ and $b$ are relatively
prime.

 In this paper, an explicit
formula for $t(\cdot |M)$ is presented for any integer matrix $M$.
In contrast to other ways, our method is based on multivariate
spline functions and is inspired by the work of Dahmen, Micchelli
and Jia \cite{dahmen1,dahmen2,Jia:discret spline,jiamagic}, who
exploited the relation between $t(\cdot|M)$ and multivariate
splines, and
 demonstrated  the piecewise structure of
 $t(\cdot|M)$.
  Moreover,  we believe
that the tool of multivariate splines, which have been developed in
this latter theory,  shed some light on  problems concerning
$t(\cdot |M )$.

The main results in this paper are  as follows. As the central
result of the paper,  an explicit formula for $t(\cdot |M)$ in terms
of multivariate splines is presented. As applications of our
formula, we first generalize the formula for partition functions in
\cite{BeckFrobenius} and give a simple proof for Popoviciu's
formula. Further, we show that the relative volume of a convex
polytope agrees with the value of the multivariate truncated power
at a certain point. We also present an efficient method for
computing the volume of convex polytopes and re-prove the volume
formula for the polytope related to empirical distributions, which
is presented in \cite{volume} as a central result.

The paper is organized as follows. To help make this paper
self-contained, we  introduce the multivariate truncated power and
the  box spline in Section 2. In Section 3, the discrete truncated
power $t(\cdot |M)$ is introduced. An explicit formula for
$t(\cdot |M)$ is described in Section 4. The proof of the formula
is presented in Section 5. In Section 6, we simplify the explicit
formula for $t(\cdot|M)$ for some special matrices $M$. The
explicit formula for partition functions in \cite{BeckFrobenius}
follows directly from our  formula. Section 7, containing two
subsections, uses multivariate truncated powers to investigate the
volume of convex polytopes and the Ehrhart quasi-polynomial.
Particularly, in Subsection 7.1, an efficient method for computing
the volume of convex polytopes is introduced, with re-proving the
volume formula for the polytope related to empirical distributions
which is the  main result in \cite{volume}; in Subsection 7.2, an
explicit formula for the Ehrhart quasi-polynomial is presented.

It is necessary here to recall some previous notations (see
\cite{jiamagic}). Throughout the paper, ${\Z}_+$ and ${\R}_+$
denote the non-negative integer and non-negative real sets
respectively. For given  sets $D_1$, $D_2$, let $1_{D_1}(D_2)=0$
if $D_1\cap D_2=\emptyset$, otherwise let $1_{D_1}(D_2)=1$.  The
linear space ${\R}^s$ is equipped with the norm $|\cdot |$ given
by $ |x|=\sum_{1\leq j\leq s}|x_j|,$ where $
x=(x_1,x_2,\ldots,x_s)\in {\R}^s$. Let $A$ and $B$ be two subsets
of ${\R}^s$. Then $A-B$ is the set of all elements in the form of
$a-b,$ where $a\in A$ and $b\in B.$ The sets $A+B$ and $cA$ are
defined analogously,  where $c\in {\R}$. The set $A\backslash B$
is the complement of $B$ in $A$. A subset $\Omega$ of ${\R}^s$ is
called a {\it cone} if $\Omega +\Omega \subseteq \Omega$ and
$c\Omega\subseteq \Omega$ for all $c>0$. If a cone is also an open
set, then we call it an {\it open cone}.  Let $Y$ be an $s\times
n$ matrix. The linear span of $Y$, denoted by ${\rm span}(Y)$, is
the set $ \{ \sum_{y\in Y}a_y y : a_y\in {\R} \mbox{ for all }
y\}. $ The cone spanned by $Y$, denoted by ${\rm cone}(Y)$, is the
set $ \{ \sum_{y\in Y}a_y y:a_y\geq 0 \mbox{ for all\ } y \}$. We
denote by ${\rm cone}^\circ (Y)$ the relative interior of ${\rm
cone}(Y).$ The {\it convex hull} of $Y$, denoted by $[Y],$ is the
set $ \{ \sum_{y\in Y}a_y y:a_y\geq 0 \mbox{ for all } y \mbox{
and } \sum_{y\in Y}a_y=1\}$. Let $M$ be an $s\times n$ matrix with
integer columns $m_1,\ldots,m_n$. We denote by $[[M]]$ the
zonotope spanned by $m_1,\ldots, m_n$, i.e., $
[[M]]:=\{\sum_{j=1}^na_jm_j: 0\leq a_j\leq 1, \forall j\}$.
Moreover, set $ [[M)):=\{\sum_{j=1}^na_jm_j: 0\leq a_j< 1, \forall
j\}$.

We  shall use  standard multiindex notations. Specifically, an
element $\alpha\in {\Z}_+^s$ is called an {\it $s$-index}, and
$|\alpha |$ is called the length of $\alpha$. Define $ z^\alpha
:=z_1^{\alpha_1}\cdots z_s^{\alpha_s}$ for $z=(z_1,\ldots,z_s)\in
{\C}^s$ and $\alpha=(\alpha_1,\ldots,\alpha_s)\in {\Z}_+^s$. Also,
set $\exp(c \alpha):=(\exp(c \alpha_1),\ldots, \exp(c \alpha_s))$,
where $c$ is a constant.
 For $z=(z_1,\ldots,z_s)\in {\C}^s$ and
$w=(w_1,\ldots,w_s)\in {(\C \setminus 0) }^s$, set
$z/w:=(z_1/w_1,\ldots,z_s/w_s)$.  We denote by
${\PP}_k:={\PP}_k({\R}^s)$ the linear space of polynomials in $s$
real variables with total degree $\leq k$, where $k\in {\Z}$. If $k$
is a negative integer, then we interpret ${\PP}_k$ as the trivial
linear space $\{ 0\}$. Moreover, we also set
${\PP}(\R^s):=\bigcup_{k=0}^\infty \PP_k(\R^s)$.
\section{Multivariate truncated powers and  box splines}
\setcounter{section}{2}

Let $M$ be an $s\times n$ real matrix with ${\rm rank}(M)=s$.
Throughout this section we always assume that the convex hull of
$M$ does not contain the origin. The multivariate truncated power
$T(\cdot |M)$ associated with $M$, introduced firstly by  Dahmen
\cite{dahemenpower},  is the distribution given by the rule
\begin{equation}\label{powerdefinition}
\int_{{\R}^s}\!\!T(x |M)\phi (x)dx\,\, =\,\, \int
_{{\R}_+^n}\!\!\phi (Mu)du , \,\,\, \phi \in {\mathscr D}({\R}^s),
\end{equation}
where ${\mathscr D}({\R}^s)$ is the space of test functions on
${\R}^s$, i.e., the space of all compactly supported and
infinitely differentiable functions on ${\R}^s$. From
(\ref{powerdefinition}), we see that the support of $T(\cdot |M)$
is ${\rm cone}(M).$ In fact, $T(\cdot |M)$ is identified with the
function (see \cite{deboorbook} pp.12)
\begin{equation}\label{Eq:truncatedvolume}
T(x|M)\,\,=\,\,\frac{{\rm vol}_{n}(M^{-1}x\cap
{\R}_+^n)}{\sqrt{|\det(MM^T)|}},
\end{equation}
where $M^{-1}x:=\{ y: My=x \}$ and ${\rm vol}_{n}(M^{-1}x\cap
{\R}_+^n)$ denotes the volume of $M^{-1}x\cap {\R}_+^n$ in
${\R}^n$. In the following, we  review some basic properties of
multivariate truncated powers. For more detailed information about
this function, the reader is  referred  to
\cite{dahemenpower,deboorbook}.

For $y=(y_1,\ldots, y_s) \in {\R}^s$ and a function $f$ defined on
${\R}^s$, we denote by $D_yf$ the directional derivative of $f$ in
the direction $y$, i.e., $ D_y=\sum_{j=1}^sy_jD_j, $ where $D_j$
denotes a partial derivative with respect to the $j$th coordinate.
The following differential formula was given in
\cite{dahemenpower}. For $y\in M$, we have $ D_yT(\cdot
|M)=T(\cdot |M\backslash y)$. More generally,
\begin{equation}\label{eq:Tre}
 D_YT(\cdot
|M)\,\,\,=\,\,\,T(\cdot |M\backslash Y)\,\,  \mbox{ for }\,\,
Y\subset M
\end{equation}
 where  $ D_Y\,:=\,\prod_{y\in Y}D_y$.

Let ${\mathscr Y}(M)$ denote the set consisting of those $Y\subset
M$ such that $M \backslash Y$ does not span ${\R}^s$ and let
$D(M)$ denote the linear space of those infinitely differentiable
complex-valued functions $f$ on ${\R}^s$ that satisfy the
following system of linear partial differential equations:
$$
 D_Yf=0,\,\, \forall\,\,\,\, Y\in
{\mathscr Y}(M).
$$
Based on the definition of $D(M)$, we can see that $D(M')\subset
D(M)$ when $M'\subset M$. It was also proved in
\cite{DahmenDM,deboor2} that $ D(M)\subseteq {\PP}_{n -s}. $

 Let $c(M)$ be the union of
all the sets cone$(M\setminus Y)$, as $Y$ runs over ${\mathscr
Y}(M)$. A connected component of ${\rm cone}^\circ (M)\setminus
c(M)$ is called a {\it fundamental $M$-cone} according to
\cite{dahmen2}.   Then $T(\cdot |M)$ agrees with some homogeneous
polynomial of degree $n-s$ in $D(M)$ on each fundamental $M$-cone.
In fact,\ $T(\cdot |M)$ is generally continuous and positive on
${\rm cone} ^\circ (M)$.

We now turn to  box splines. The  box spline $B(\cdot |M)$
associated with $M$ is the distribution given by the rule
\cite{deboor1,deboor2}
\begin{equation}\label{Eq:boxspline}
\int_{{\R}^s}\!\!B(x |M)\phi (x)dx \,\,\,=\,\,\, \int
_{[0,1)^n}\!\!\phi (Mu)du ,\,\,\,\, \phi \in {\mathscr D}({\R}^s).
\end{equation}
According to (\ref{Eq:boxspline}), the support of $B(\cdot|M)$ is
$[[M]]$, and hence, we have (see \cite{deboorbook}, pp.33)
$$
\{ j\in {\Z}^s: B(x-j|M)\not= 0\}\,\,=\,\, {\Z}^s\cap (x-[[M))).
$$
By taking $\phi={\exp}(-i\xi\cdot )$ in (\ref{Eq:boxspline}), we
obtain the Fourier transform of $B(\cdot |M)$ as
$$
\widehat{B}(\zeta
|M)\,\,=\,\,\prod_{j=1}^n\frac{1-\exp(-i\zeta^Tm_j)}{i\zeta^Tm_j},\,\,\,\,
\zeta \in {\C}^s.
$$
For more detailed information about box splines, the reader is
referred to \cite{deboorbook}.

\section{Discrete truncated powers}

\setcounter{section}{3}

Let $M$ be an $s\times n$ matrix with integer columns
$m_1,\ldots,m_n$ and suppose that $[M]$ does not contain the
origin. From the definition of $t(\cdot |M)$ given in Section 1,
we have
\begin{equation}\label{eq:discre}
\sum_{\alpha \in {\Z}^s}\varphi (\alpha )t(\alpha
|M)\,\,=\,\,\sum_{\beta\in {\Z}_+^n}\varphi (M\beta),\,\,\, \mbox{
for any }\varphi \in {\mathscr D}({\R}^s).
\end{equation}
A comparison between (\ref{eq:discre}) and (\ref{powerdefinition})
shows that (\ref{eq:discre}) is a discrete version of
(\ref{powerdefinition}). This observation was the motivation to
designate $t(\cdot |M)$ a {\it discrete truncated power}. Also the
identity
\begin{equation}\label{Formula:three functions closely related}
T(x|M)\,\,=\,\,\sum_{\alpha\in {\Z}^s}t(\alpha |M)B(x-\alpha | M)
\end{equation}
established in \cite{dahmen1} shows that these three functions,
$T(\cdot |M)$,  $B(\cdot|M)$ and $t(\cdot |M)$, are closely
related.

We denote by $S$ the linear space of all complex functions on
${\Z}^s$. Given $y\in {\Z}^s$, the backward difference operator
$\nabla_y$ is defined by the rule $ \nabla_yf:=f-f(\cdot -y),$
where $f\in S. $ More generally, for an integer matrix $Y$, we
define
$$ \nabla_Y\,\,:=\,\,\prod_{y\in Y}\nabla_y.$$
In \cite{dahmen2}, the
following difference formula was given:
$$
\nabla_yt(\cdot|M)\,\,=\,\,t(\cdot |M\setminus y), \mbox{ for any
} y\in M.
$$
More generally,
$$ \nabla_Yt(\cdot |M)\,\,=\,\,t(\cdot |M\setminus
Y),\mbox{ for any } Y\subseteq M.
$$
 Given $\theta=(\theta_1,\ldots, \theta_s)\in
({\C}\setminus 0)^s,$ we  set
$$ M_{\theta}\,\,:=\,\,\{ y\in
M: \theta^y=1\}.
$$
 Let
 $$ A(M)\,\,:=\,\,\{ \theta\in ({\C}\setminus
0)^s : {\rm span}(M_{\theta})={\R}^s \}.$$
 As  pointed out in
\cite{DaMi on the solution}, $\theta\in A(M)$ if and only if it
has the form
\begin{equation}\label{the term of A(M)}
\theta\,\, =\,\,\exp (2\pi i \alpha/\det Y),
\end{equation}
for some $Y\in {\mathscr B}(M)$,  and the vector $\alpha\in \Z^s$
satisfying  $ Y^T\alpha=|\det Y|L$,  where $L\in \Z^s \cap[[Y^T))$,
i.e., is a lattice point in the parallelepiped determined by $Y^T$.
Here,
$${\mathscr B}(M)=\{ Y\subseteq M:\#Y=s,\,\,
{\rm span}(Y)={\R}^s\}.
$$
We let
$$
\nabla(M)\,\,:=\,\,\{f\in S : \nabla_Y f=0,\,\,\, \mbox{for all
}\,\,\, Y\in {\mathscr Y}(M)\}.
$$
In \cite{DaMi on the solution}, Dahmen and Micchelli showed that a
sequence $f\in \nabla (M)$ if and only if it has the form $
f(\alpha)=\sum_{\theta\in
A(M)}\theta^{\alpha}p_{\theta}(\alpha),\alpha\in {\Z}^s, $ where
$p_{\theta}\in D(M_{\theta})$ for each $\theta\in A(M)$. For a
fundamental $M$-cone $\Omega$, we set
$$
 v(\Omega |M)\,\,:=\,\,{\Z}^s \cap
(\Omega-[[M]]).
$$ Dahmen and Micchelli also proved the following
 result in \cite{dahmen2}:
\begin{theorem}{\rm (\cite[Thm.3.1]{dahmen2})}\label{Theorem:DaMi,truncated power structure}
 Let $M=\{m_1,\ldots,m_n\}$ be an $s\times n$ integer matrix and suppose
  that $M$ spans ${\R}^s$ and the convex hull of $M$
does not contain the origin. Then for any fundamental $M$-cone
$\Omega,$ there exists a unique element $f_{\Omega}(\cdot |M)\in
\nabla (M)$ such that $f_{\Omega}(\cdot |M)$ agrees with $t(\cdot
|M)$ on $v(\Omega |M).$ Moreover, $f_{\Omega}(\cdot |M)$ has the
following properties: for any $x\in \Omega$ such that
$$
v(x|M)\cap {\rm cone}(M)=\{ 0\},
$$
$f_{\Omega}(\cdot |M)$ is uniquely determined by
$$
f_\Omega(\alpha|M)=\delta_{0\alpha },\,\, \alpha \in v(x|M),
$$
and satisfies the relation
\begin{equation}\label{Equation:the relation}
f_{\Omega}(\alpha |M)=(-1)^{n-s}f_{\Omega}(-\alpha
-\sum_{j=1}^nm_j|M),\,\,\,\, \alpha \in {\Z}^s.
\end{equation}
\end{theorem}

 We denote by $\theta^{()}p_{\theta }$ the
complex sequence given by $\alpha \mapsto \theta^\alpha p_{\theta
}(\alpha ),$ with $\alpha \in {\Z}^s$, where $p_{\theta}\in
{\PP}({\R}^s)$ and $\theta\in ({\C}\setminus 0)^s$.
 Let $E$ denote the linear space
of all sequences of the form $f=\sum_{\theta}\theta^{()}p_{\theta}.$
It is easily seen that each $f\in E$ can be written uniquely in the
form $\sum_{\theta}\theta^{()}p_{\theta}$. Hence, we can define
mappings $P_{\theta}$ from $E$ into ${\PP}({\R}^s)$ as $P_{\theta}:
f \mapsto p_{\theta}$.

Based on Theorem \ref{Theorem:DaMi,truncated power structure},
$f_{\Omega}(\alpha |M)\in \nabla(M)$. So, it can be written in the
form of
\begin{equation}\label{eq:fomega}
f_{\Omega}(\alpha |M)\,\,=\sum_{\theta\in
A(M)}\theta^{\alpha}p_{\theta}(\alpha),\,\,\,\,\,\,\, \alpha\in
{\Z}^s,
\end{equation}
where $p_{\theta}(\alpha )\in D(M_{\theta})$. According to the
definition of $P_{\theta}$, one has $P_{\theta }f_{\Omega }(\alpha
|M)=p_{\theta}(\alpha),$ for any $\theta\in A(M).$ Let
$e:=(1,1,\ldots,1)\in {\Z}^s. $ As, for any $y\in M$, $e^y=1,$ we
see that $e\in A(M)$.  In fact, $P_ef_{\Omega}(\cdot|M)$ is the
polynomial part of $f_{\Omega}(\cdot|M).$ In \cite{dahmen2}, Dahmen
and Micchelli showed that the following equation holds: for $x\in
\Omega,$
\begin{equation}\label{damiequation}
T(x|M)\,\,=\,\,P_ef_{\Omega}(x|M)\,\,+\sum_{0<|u|\leq
n-s}D^uP_ef_{\Omega}(x|M)(-i)^{|u|}D^u\widehat{B}(0|M)/u!.
\end{equation}
From this equation, the following theorem was proved in
\cite{dahmen2}.
\begin{theorem}{\rm (\cite[Prop.5.3]{dahmen2})}\label{Theorem:the leading terms}
Under the conditions of Theorem \ref{Theorem:DaMi,truncated power
structure}, the leading homogeneous term of $P_{e}f_{\Omega}(\cdot
|M)$ agrees on $\Omega$ with $T(\cdot |M).$
\end{theorem}

\section{An explicit formula for discrete truncated powers}

\setcounter{section}{4}

The objective of this section is to present an explicit formula for
discrete truncated powers $t(\cdot |M)$ which is the central result
of the paper. Throughout the rest of the paper, we use $\Omega$ to
denote any particular fundamental $M$-cone. According to Theorem
\ref{Theorem:DaMi,truncated power structure}, the discrete truncated
power $t(\cdot|M)$ agrees, on $v(\Omega |M)$,  with an element in
$\nabla (M)$, which is denoted as $f_\Omega(\cdot |M)$. We use
$\overline{\Omega}$ to denote the closure of $\Omega$. From Lemma
4.4 in  \cite{jiamagic}, we have $\overline{\Omega}\cap \Z^s \subset
v(\Omega|M)$. Note that $t(\alpha |M)=0,$ when $\alpha\in
{\Z}^s\setminus \overline{{\rm cone}(M)}$. Hence, to present an
explicit formula for $t(\cdot |M),$ we only need an explicit formula
for $f_\Omega(\cdot |M)$.

In the following theorem, we  give the polynomial part of
$f_{\Omega}(\cdot|M)$, i.e., $P_ef_{\Omega}(\cdot|M)$.
\begin{theorem}\label{Th:FormP}
Suppose that $p_{\mu,\Omega}$ is the homogeneous polynomial part of
degree $n-s-\mu$   of
 $P_ef_{\Omega}(x|M)$.
Under the conditions of Theorem \ref{Theorem:DaMi,truncated power
structure}, when $x\in \Omega$, we have
\begin{eqnarray}\label{eq:thre}
p_{\mu,\Omega}(x)=
\begin{cases}
T(x|M) ,\hspace{5cm} \mu=0,\\
-\sum_{j=0}^{\mu-1}\left(\sum_{|u|=k-j}D^vp_{j,\Omega}(x)
(-i)^{|u|}D^u\widehat{B}(0|M)/u!\right),
\end{cases}\\
1\leq \mu \leq n-s.\nonumber
\end{eqnarray}
\end{theorem}
\begin{proof}
 Rewriting of Equation (\ref{damiequation}) yields
\begin{equation}\label{eq:dami}
\begin{split}
&T(x|M)\,\,-\sum_{0<|u|\leq
n-s}D^uP_ef_{\Omega}(x|M)(-i)^{|u|}D^u\widehat{B}(0|M)/u! \\
=&\,\,P_ef_{\Omega}(x|M),\,\,\, x\in \Omega .
\end{split}
\end{equation}
According to (\ref{eq:dami}), $p_{0,\Omega}(x)$ agrees with
$T(x|M)$ on $\Omega$. Since both sides of (\ref{eq:dami}) are
polynomials, we can rewrite (\ref{eq:dami}) as
\begin{equation}\label{eq:prmid}
p_{0,\Omega}(x)-\sum_{0<|u|\leq
n-s}D^u\sum_{\mu=0}^{n-s}p_{\mu,\Omega}(x)(-i)^{|u|}D^u\widehat{B}(0|M)/u!
=\sum_{\mu=0}^{n-s}p_{\mu,\Omega}(x),\,\,\, x\in \R^s.
\end{equation}
A comparison between both sides of the equation above shows that
$$p_{\mu,\Omega}(x)\,\,=\,\, -\sum_{j=0}^{\mu-1}\left(\sum_{|u|=\mu-j}D^up_{j,\Omega}(x)
(-i)^{|u|}D^u\widehat{B}(0|M)/u!\right),\,\, 1\leq \mu\leq n-s.$$
\end{proof}
Based on (\ref{eq:fomega}), to write down an explicit formula for
$t(\cdot |M)$ on $\Omega$, we only need  a formula for
$P_{\vartheta} f_{\Omega}(\cdot|M)$ for each fixed $\vartheta\in
A(M)$. To describe the formula conveniently, we suppose
$M\setminus
 M_{\vartheta}=(m_1,\ldots,m_\kappa)$, where $\kappa =\#M-\#M_{\vartheta}$.
 We select $r$ as the least integer in the set
$ \{r: (\vartheta^m)^r=1 \mbox{ for all } m \in M\setminus
M_{\vartheta}\}$. (According to (\ref{the term of A(M)}), the set
is non-empty, since there is at least an integer,  $\prod_{Y\in
{\mathscr B}(M)}|\det(Y)|$, in it.)
 We  use $q_{\mu,r}^{\vartheta}(x)$, $\mu\in {\Z}_+$, to denote homogeneous
polynomials satisfying  the following condition:
\begin{equation}\label{eq:qform}
q_{ \mu ,r}^{\vartheta}(x)=\sum_{ j_1+\cdots+j_\kappa =\mu
}\frac{1}{r^\kappa}D_{m_1}^{j_1}\cdots
D_{m_\kappa}^{j_\kappa}T(x|M_{{\vartheta}})\prod_{i=1}^\kappa
\frac{s_{j_i+1}({\vartheta}^{-m_i})}{(j_i+1)!},\,\,\,\, x\in
\Omega.
\end{equation}
Throughout the rest of the paper, we set
$$s_j(x)\,\,:=\,\,(-1)^j(x+2^jx^2+\cdots+(r-1)^jx^{r-1})$$ and
$\widetilde{M}_{r}:=\{\widetilde{m}_1,\ldots,\widetilde{m}_n\},$
where if ${m_i}\in M_{\vartheta}$, then $\widetilde{m}_i=m_i,$
otherwise, $\widetilde{m}_i=rm_i$.
\begin{theorem}\label{Th:FormPtheta}
Suppose that $p^{\vartheta}_{\mu,\Omega}$ is the homogeneous
polynomial part of degree $n-s-\kappa-\mu$ of
 $P_{\vartheta}f_{\Omega}(\cdot|M)$.
 Under the conditions
of Theorem \ref{Theorem:DaMi,truncated power structure}, for each
fixed $\vartheta\in A(M)$, we have
\begin{equation*}
\begin{split}
p^{\vartheta}_{0,\Omega}(\cdot)&=q_{0,r}^{\vartheta}(\cdot),\\
p^{\vartheta}_{\mu,\Omega}(\cdot)&=q_{\mu,r}^{\vartheta}(\cdot)
-\sum_{j=0}^{\mu-1}\left(\sum_{|u|=\mu-j}D^vp^{\vartheta}_{j,\Omega}(\cdot)
(-i)^{|u|}D^v\widehat{B}(0|\widetilde{M}_{r})/u!\right),\\
&\hspace{7cm} 1\leq \mu\leq n-s-\kappa,
\end{split}
\end{equation*}
where   $r$ is the least positive integer such that
$(\vartheta^m)^r=1$ holds for all $m\in M\setminus M_{\vartheta}$.
\end{theorem}

Combining Theorem \ref{Th:FormPtheta} and
$q_{0,r}^{\vartheta}(x)=T(x |M_{{\vartheta}})\prod_{w\in M\setminus
M_{{\vartheta}}}\frac{1}{1-{\vartheta}^{-w}}$,  $x\in \Omega$, we
can easily extend Theorem \ref{Theorem:the leading terms}  as
follows:
\begin{theorem}{\rm (\cite{wangxu})}\label{Theorem:leadingpart of all }
Under the conditions of Theorem \ref{Theorem:DaMi,truncated power
structure}, the leading part of $P_{{\vartheta}}f_{\Omega}(\cdot
|M)$ agrees with
$$T(\cdot
|M_{{\vartheta}}) \prod_{w\in M\setminus
M_{{\vartheta}}}\frac{1}{1-{\vartheta}^{-w}}$$
 on $\Omega $.
\end{theorem}

An alternative method for proving Theorem \ref{Theorem:leadingpart
of all } is also given in \cite{wangxu}.
\begin{remark}
 Since $P_{\vartheta}f_{\Omega}(\cdot |M)$ is a polynomial and
$\overline{\Omega}\cap \Z^s \subset v(\Omega|M)$ (see Lemma 4.4 in
\cite{jiamagic}), Theorem \ref{Theorem:the leading terms} and
Theorem \ref{Theorem:leadingpart of all } also hold on $\overline
{\Omega}$.
\end{remark}

\begin{remark}
In \cite{brion},  Brion and Vergne  give formulas for the lattice
point enumerator of a convex rational polytope in terms of certain
Todd differential operators. These are interesting connections to
topology. The salient difference in approaches lies in our explicit
use of splines.
 Although  the formula presented in \cite{brion}  might
eventually be equivalent our, the analysis of the equivalence
between our different formulas would be the content of another
paper.
\end{remark}

\begin{example} We consider the number of non-negative integer
solutions of the linear equations
 $$
\begin{pmatrix}
 3 & 2 & 1 & 0 \\
 0 & 1 & 2 & 2
\end{pmatrix}\begin{pmatrix}
x_1 \\
x_2 \\
x_3\\
x_4
\end{pmatrix} = \begin{pmatrix}
k_1 \\
k_2
\end{pmatrix}.
$$
Set $
 M=\begin{pmatrix}
 3 & 2 & 1 & 0 \\
 0 & 1 & 2 & 2
\end{pmatrix}$ and ${\bf k}=\begin{pmatrix}
k_1 \\
k_2
\end{pmatrix}$.
A simple computation shows that $A(M)=\{(1,1)\}\cup A_3(M)\cup
A_2(M)$, where
\begin{eqnarray*}
 A_3(M)&=&\{ (1,-1),\, (\exp{(2\pi i/3)},\,\exp{(2\pi i/3)}),\,(\exp{(4\pi i/3)},\,\exp{(4\pi i/3)})\},\\
 A_2(M)&=&\{(-1,1),\,(\exp{(2\pi i/3)},1), (\exp{(2\pi i/3)},-1), (\exp{(4\pi
 i/3)},1),\\
 & & (\exp{(4\pi i/3)},-1),
 (i,-1),\, (\exp{(3\pi i/2) },-1),\,(\exp{(2\pi i/3)},\exp{(5\pi
 i/3)}),\\
& &(\exp{(4\pi i/3)},\exp{(\pi i/3)})\}.
\end{eqnarray*}
For $\theta\in A_3(M)$, we have $\#M_{\theta}=3$ while $\#
M_{\theta}=2$ provided that $\theta\in A_2(M)$. For the matrix $M
$, there are $3$ fundamental $M$-cones, i.e., $\Omega_1:={\rm
cone}^\circ(M_{12})$, $\Omega_2:={\rm cone}^\circ(M_{23})$ and
$\Omega_3:={\rm cone}^\circ(M_{34})$, where
$M_{12}=\begin{pmatrix}
3 &2\\
0 & 1
\end{pmatrix},
M_{23}=\begin{pmatrix}
2 &1\\
1 & 2
\end{pmatrix}\,\, {\rm and }\,\,
M_{34}=\begin{pmatrix}
1 &0\\
2 & 2
\end{pmatrix}.$

Using Theorem \ref{Th:FormPtheta}, we  obtain the explicit formula
for $t(\cdot|M)$ as follows:  when ${\bf k}\in v(\Omega_1,M)$,
\begin{eqnarray*}
& &t({\bf k}|M)=\\
& & {k_2^2}/{24}+ {5k_2}/{24}+ {11}/{48}+{5}/{48}+
\left({(3+\sqrt{3}i)}k_2/{36}
+{(17+3\sqrt{3}i)}/{72}\right)\\
& &\cdot \exp({(2\pi i k_1+2\pi i k_2)/3})
+(-1)^{k_2}({k_2}/{24}+5/48)+
 \exp({(4\pi i k_1+4\pi i
k_2)/3})\\
& &
\left({(3-\sqrt{3}i)}k_2/{36}+{(17-3\sqrt{3}i)}/{72}\right)+{\exp({2\pi
i k_1/3 })}/{18} -{\sqrt{3}}/{18}i\exp({2\pi i
k_1/3})(-1)^{k_2}\\
& & +\exp({4\pi i k_1/3})/{18}+{\sqrt{3}}/{18} i
(-1)^{k_2}\exp({4\pi i k_1/3}) +\exp({(2\pi ik_1+5\pi i
k_2)/3})({1}/{24}\\
& &+{\sqrt{3}}/{72}i)+\exp({(4\pi i k_1+\pi i
k_2)/3})({1}/{24}-{\sqrt{3}}/{72}i);
\end{eqnarray*}

when ${\bf k}\in v(\Omega_2,M)$,
\begin{eqnarray*}
& &t({\bf k}|M)=\\
 & &  {k_1k_2}/{18}-{(k_1^2+k_2^2)}/{72}
+{7k_2}/{72}+{k_1}/{18}+{23}/{108}
+(-1)^{k_2}({k_2}/{24}+{5}/{48})\\
& &+ \exp({(2\pi i k_1+2\pi i
k_2)/3})\left({(3+\sqrt{3}i)}/{108}\cdot {(2k_1-k_2)}+{19}/{216}+{\sqrt{3}}/{24}i\right)\\
& &+
 \exp({(4\pi i k_1+4\pi i k_2)/3})\left({(3-\sqrt{3}i)}/{108}\cdot {(2k_1-k_2)}+{19}/{216}-{\sqrt{3}}/{24}i\right)\\
& &+ {(-1)^{k_1}}/{16}+ {1}/{18}\exp({2\pi i k_1/3 })
-{\sqrt{3}}/{18}i\exp({2\pi i
k_1/3})(-1)^{k_2} +\exp({4\pi i k_1/3})/18\\
& &+{\sqrt{3}}/{18}i(-1)^{k_2}\exp({4\pi i k_1/3 })
+{1}/{8}\exp({\pi i k_1/2
})(-1)^{k_2}+{1}/{8}(-1)^{k_2}\exp({3\pi i k_1/2})\\
 & &+
({1}/{24}+{\sqrt{3i}}/{72})\exp({(2\pi i k_1+5\pi i
k_2)/3})+({1}/{24}-{\sqrt{3}i}/{72})\exp({(4\pi i k_1+\pi i
k_2)/3});
\end{eqnarray*}

when ${\bf k}\in v(\Omega_3,M)$,
\begin{eqnarray*}
& &t({\bf k}|M)=\\
& &{k_1^2}/{24}+{k_1}/{4}+{47}/{144}
+(-1)^{k_2}({k_1}/{12}+{1}/{4})+(-1)^{k_1}/16+\exp({2\pi
i k_1/3})/18\\
& &-({\sqrt{3}}/{18})i \exp({2\pi i k_1/3})(-1)^{k_2}
+({1}/{18})\exp({4\pi i k_1/3})+ ({\sqrt{3}}/{18})i\exp({4\pi
i k_1/3 })(-1)^{k_2}\\
 & &+({1}/{8})\exp({\pi i k_1/2 })(-1)^{k_2}
+({1}/{8})\exp({3\pi i k_1/2})(-1)^{k_2}.
\end{eqnarray*}
\end{example}

\section{Proof of Theorem \ref{Th:FormPtheta}}
\setcounter{section}{5}

In this section, we shall prove Theorem \ref{Th:FormPtheta}. The
proof begins with several lemmas, designed to make the proof more
readable. The main idea for proving Theorem \ref{Th:FormPtheta} is
to generalize  Equation (\ref{damiequation}), with the equation
playing an important role in obtaining the polynomial part for
$f_{\Omega}(\cdot |M)$. This is the plan. Our discussion will be
broken into three steps.

First, we introduce a set. Let
$$c(\Omega,H)\,\,:=\,\,\bigcap_{h\in H}\bigcap_{b_h\in [0,1]
}(\Omega+b_hh),$$
 where $H$
is a finite set of real $s$-vectors. For example, if we set
$H:=\{1\}$ and $\Omega:=(0,+\infty)$, then
$c(\Omega,H)=(1,+\infty)$. We have
\begin{lemma}\label{Le:asetpro} Suppose that $H$
is a finite set of real $s$-vectors. Then for the fundamental
$M$-cone $\Omega$, the set $c(\Omega,H)\subset \Omega$ and the
volume of $c(\Omega,H)$ in ${\R}^s$ is infinity.
\end{lemma}
\begin{proof}
According to the definition of $c(\Omega,H)$, one has
$$c(\Omega,H)=\bigcap_{h\in H}\bigcap_{b_h\in [0,1]
}(\Omega+b_hh)\subset \bigcap_{h\in H,b_h=0}(\Omega+b_hh)=\Omega.$$
Based on the definition of the fundamental $M$-cone, there exist
real $s$-vectors ${g}_1,\ldots,{ g}_\omega$, $\omega \in {\N} $,
such that ${\rm span}\{{ g}_1,\ldots,{ g}_\omega\}={\R}^s$ and
$\Omega=\{a_1{ g}_1+\cdots+a_\omega{ g}_\omega:a_i>0\}$. Hence, for
any $h\in H$, there exist
$\widetilde{a}_1(h),\ldots,\widetilde{a}_\omega(h)\in {\R}$ such
that $h=\widetilde{a}_1(h){ g}_1+\cdots+\widetilde{a}_\omega (h){
g}_\omega $. Then
\begin{eqnarray*}
c(\Omega,H)=\bigcap_{h\in H}\bigcap_{ b_h\in[0,1]}
\bigcup_{a_1,\ldots,a_\omega>0}\{\sum_{i=1}^\omega
(a_i+b_h\widetilde{a}_i(h)){\bf g}_i\}.
\end{eqnarray*}

 Let $a(h):=\max \{
|\widetilde{a}_1(h)|,\ldots, |\widetilde{a}_\omega (h)|\}$ and let
$a(H):=\max\{ a(h): h\in H\}$. A simple observation is that $a_i+b_h
\widetilde{a}_i(h)>0$ for any $b_h\in [0,1]$, $h\in H$, provided
that $a_i>a(H)$. Hence, we have
$${\bf
S}:= \{{ x}: { x}=a_1h_1+\cdots+a_\omega h_\omega,a_i\in {\R},
a_i>a(H)\} \subset c(\Omega,H).$$ Note that the volume of ${\bf S}$
in ${\R}^s$ is infinite.
 So,
the volume of $c(\Omega,H)$ in ${\R}^s$ is also infinite.
\end{proof}

Second, we shall generalize Equation (\ref{damiequation}), since
it is of importance in obtaining the explicit formula for
$P_ef_\Omega(x|M)$. We recall some notations and definitions.  Let
$\vartheta$ be an arbitrary but a fixed vector in $A(M)$. Recall
that $r$ is the least positive integer such that
$(\vartheta^m)^r=1$ holds for all $m\in M\setminus M_{\vartheta}$.
Let
$\widetilde{M}_{r}=\{\widetilde{m}_1,\ldots,\widetilde{m}_n\},$
where if ${m_i}\in M_{\vartheta},$ then $\widetilde{m}_i=m_i,$
otherwise, $\widetilde{m}_i=rm_i$.  Before extending Equation
(\ref{damiequation}), we first generalize (\ref{Formula:three
functions closely related}) as follows:

\begin{prop}\label{le:tBT}
\begin{equation}\label{eq:le}
 \sum_{\alpha\in {\Z}^s}\vartheta^{-\alpha}t(\alpha
|M)B(x-\alpha |\widetilde{M}_{r})\,\,=\! \sum_{0\leq
r_1,\ldots,r_\kappa<r}\vartheta^{-\sum_{i=1}^\kappa
m_ir_i}T(x-\sum_{i=1}^\kappa m_ir_i|\widetilde{M}_{r}).
\end{equation}
\end{prop}
\begin{proof} Using the fact that
$t(\alpha|M)=\#\{\beta\in\Z^n_+:M\beta=\alpha\}$, we see that
\begin{eqnarray}\label{TrunctedBSpline}
& &\sum_{\alpha\in {\Z}^s}{\vartheta}^{-\alpha}t(\alpha
|M)B(x-\alpha |\widetilde{M}_{r}) =\sum_{\beta\in
{\Z}_+^n}{\vartheta}^{-M\beta}B(x-M\beta|\widetilde{M}_{r})= \nonumber\\
& & \sum_{\beta_1,\ldots,\beta_{n}\in {\Z}_+
}{\vartheta}^{-\sum_{i=1}^\kappa
m_i\beta_i}B(x-\sum_{i=1}^nm_i\beta_i|\widetilde{M}_{r}).
\end{eqnarray}
We can write  $\beta_i$ as $\gamma_i r+r_i$, where $0\leq r_i<r$,
$\gamma_i, r_i\in {\Z}_+$ for $i\leq \kappa$.  Substituting
$\beta_i=\gamma_i r+r_i$, where $i\leq \kappa$, into
(\ref{TrunctedBSpline}) and noting that ${\vartheta}^{rm_i}=1,$ we
obtain that
\begin{equation*}
\begin{split}
 &\sum_{0\leq r_1,\ldots,r_\kappa<r}\,
\sum_{\gamma_1,\ldots,\gamma_\kappa\in {\Z}_+} \,
\sum_{\beta_{\kappa+1},\ldots,\beta_n\in
{\Z}_+}{\vartheta}^{-\sum_{i=1}^\kappa m_ir_i}
B(x-\sum_{i=1}^\kappa m_ir_i\\
&-r\sum_{i=1}^\kappa m_i\gamma_i -\sum_{i=\kappa+1}^n
m_{i}\beta_{i} |\widetilde{M}_{r})\\
 &=\sum_{0\leq
r_1,\ldots,r_\kappa<r}{\vartheta}^{-\sum_{i=1}^\kappa m_ir_i}
\cdot  \sum_{\beta\in
{\Z}_+^n}B(x-\widetilde{M}_r\beta-\sum_{i=1}^\kappa
m_ir_i|\widetilde{M}_{r})\\
& =\sum_{0\leq
r_1,\ldots,r_\kappa<r}{\vartheta}^{-\sum_{i=1}^\kappa
m_ir_i}T(x-\sum_{i=1}^\kappa m_ir_i|\widetilde{M}_{r}).
\end{split}
\end{equation*}

The last equation follows from $\sum_{\beta\in
{\Z}_+^n}B(x-\widetilde{M}_r\beta|\widetilde{M}_r)=T(x|\widetilde{M}_r)$,
which can be obtained directly from the definitions of $B(\cdot|M)$
and $T(\cdot |M)$.

 Now we arrive at
$$
\sum_{\alpha\in {\Z}^s}{\vartheta}^{-\alpha}t(\alpha |M)B(x-\alpha
|\widetilde{M}_{r})=\sum_{0\leq
r_1,\ldots,r_\kappa<r}{\vartheta}^{-\sum_{i=1}^\kappa
m_ir_i}T(x-\sum_{i=1}^\kappa m_ir_i|\widetilde{M}_{r}).
$$ \end{proof}

In fact, if we set $\vartheta=e$, then (\ref{eq:le}) is reduced to
(\ref{Formula:three functions closely related}). So, as said
above, (\ref{eq:le}) can be considered as a generalization of
(\ref{Formula:three functions closely related}). To simplify the
term $\vartheta^{-\alpha}t(\alpha |M)B(x-\alpha
|\widetilde{M}_{r})$ in (\ref{eq:le}), we need to study
$\theta/\vartheta$ where $\theta\in A(M)$.
\begin{prop}\label{le:theta}
{\it For any ${\theta}\in A({M}),$ we have ${\theta}/\vartheta\in
A(\widetilde{M}_r)$ and $D(M_{{\theta}})\subset
D((\widetilde{M}_r)_{\theta/\vartheta})$.}
\end{prop}
\begin{proof} The definition of $\widetilde{M}_r$ shows that
${\vartheta}^m={\vartheta}^{-m}=1$ for any $m\in \widetilde{M}_r.$
Hence, $(\widetilde{M}_r)_{{\theta}}=
(\widetilde{M}_r)_{{\theta}/\vartheta}$ for any $\theta\in A(M)$.
Also, $rM_{\theta}\subset (\widetilde{M}_r)_\theta$ since
$\theta^m=1$ implies $\theta^{rm}=1$.
 So, we have ${\rm
span}((\widetilde{M}_r)_{{\theta}/\vartheta})={\rm
span}((\widetilde{M}_r)_{{\theta}})={\rm span}({M}_{{\theta}})
={\R}^s$ which implies ${\theta}/\vartheta\in A(\widetilde{M}_r)$.
To prove that $D(M_{{\theta}})\subset D((\widetilde{M}_r)_{{\theta}/
\vartheta})$, we only need to show that $D(M_{{\theta}})\subset
D((\widetilde{M}_r)_{{\theta}})$ since
$(\widetilde{M}_r)_{{\theta}}=
(\widetilde{M}_r)_{{\theta}/\vartheta}$. Select an $f\in
D(M_{\theta})$ and consider $D_{\widetilde{Y}} f$, where
$\widetilde{Y}=(\tilde{y}_1,\ldots,\tilde{y}_w) \in {\mathscr
Y}((\widetilde{M}_r)_{\theta}) $. We claim that
$D_{\widetilde{Y}}f=0$. We set $Y=(y_1,\ldots,y_w)$, where
$y_j=\tilde{y}_j$ if $\tilde{y}_j\in M_{\vartheta} $ otherwise
$y_j=\tilde{y}_j/r$. We can see that $Y\in {\mathscr Y}(M_\theta)$
and hence $D_{\widetilde Y}f=c_0D_Yf=0$, where $c_0$ is a non-zero
constant. So, we have $f\in D((\widetilde{M}_r)_{{\theta}})$, which
implies that $D(M_{{\theta}})\subset
D((\widetilde{M}_r)_{{\theta}})$.
\end{proof}

We now generalize Equation (\ref{damiequation}).
 Let $H_{{\vartheta}}:=\{
(r-1)m_i:1\leq i\leq \kappa \}$ and
$$Q_{{\vartheta},r}(x):=\sum\limits_{0\leq
r_1,\ldots,r_\kappa<r}{\vartheta}^{-\sum_{i=1}^\kappa
m_ir_i}T(x-\sum_{i=1}^\kappa m_ir_i|\widetilde{M}_{r}).$$ Then we
have the following result.
\begin{lemma}\label{le:4}
When $x\in c(\Omega,H_{{\vartheta}})$,
\begin{equation}\label{Eq:relationproof}
Q_{{\vartheta},r}(x)=P_{{\vartheta}}f_{\Omega}(x|M)+\sum_{|u|=1}^{n-\kappa-s}D^uP_{{\vartheta}}
f_{\Omega}(x|M)(-i)^{|u|}D^u\widehat{B}(\cdot|\widetilde{M}_{r})(0)/u!.
\end{equation}
\end{lemma}
\begin{proof}
By Proposition \ref{le:tBT} we have
$$
Q_{{\vartheta},r}(x)\,\,=\,\, \sum_{\alpha\in
{\Z}^s}\vartheta^{-\alpha}t(\alpha |M)B(x-\alpha
|\widetilde{M}_{r}).
$$
So, to this end,  we only need to prove that
\begin{equation*}\label{eq:prlema541}
\sum_{\alpha\in {\Z}^s}\vartheta^{-\alpha}t(\alpha |M)B(x-\alpha
|\widetilde{M}_{r})=P_{{\vartheta}}f_{\Omega}(x|M)+\!
\sum_{|u|=1}^{n-\kappa-s}
D^uP_{{\vartheta}}f_{\Omega}(x|M)(-i)^{|u|}D^u\widehat{B}(\cdot|\widetilde{M}_{r})(0)/u!.
\end{equation*}
  The definition of $c(\Omega, H_{{\vartheta} })$ shows that
$v(x|\widetilde{M}_r)\subset v(\Omega |M)$ provided  that $x\in
c(\Omega, H_{{\vartheta} }).$
 Noting that
$\{ j\in {\Z}^s : B(x-j|\widetilde{M}_r)\not=
0\}=v(x|\widetilde{M}_r)$ (see \cite{deboorbook}, pp.33), we have
\begin{equation}\label{eq:peq1}
\sum_{\alpha\in v(\Omega|M)}{\vartheta}^{-\alpha}t(\alpha
|M)B(x-\alpha |\widetilde{M}_{r})\,\,=\,\, \sum_{\alpha \in
{\Z}^s}{\vartheta}^{-\alpha}t(\alpha |M)B(x-\alpha
|\widetilde{M}_{r}),
\end{equation}
where $x\in c(\Omega, H_{{\vartheta} }).$
 Hence,
\begin{equation}\label{eq:le541}
\begin{split}
&\sum_{\alpha\in v(\Omega|M)}\!\! {\vartheta}^{-\alpha}t(\alpha
|M)B(x-\alpha |\widetilde{M}_{r})= \sum_{\alpha \in
{\Z}^s}{\vartheta}^{-\alpha}t(\alpha
|M)B(x-\alpha |\widetilde{M}_{r})\\
 &=\sum_{\alpha \in {\Z}^s}
\sum_{\theta\in A(M)\setminus \vartheta}(\theta/\vartheta)^\alpha
P_{{\theta}}f_{\Omega}(\alpha |M)B(x-\alpha |\widetilde{M}_{r})+
\sum_{\alpha \in {\Z}^s}P_{{\vartheta}}f_{\Omega}(\alpha
|M)B(x-\alpha
|\widetilde{M}_{r})\\
 &= \sum_{\alpha \in {\Z}^s}P_{{\vartheta}}f_{\Omega}(\alpha
|M)B(x-\alpha |\widetilde{M}_{r})\,\,\,\,\,\,\,\,\,\, {\rm
where}\,\,\, x\in c(\Omega, H_{{\vartheta} }).
\end{split}
\end{equation}
The last equation follows from $\sum_{\alpha \in
{\Z}^s}\rho(\alpha)B(x-\alpha |\widetilde{M}_r)\equiv 0$ for any
$\rho(\alpha) \in E(\widetilde{M}_r)$ (see Proposition 5.2 in
\cite{dahmen2}) and $\sum_{\theta\in A(M)\setminus
\vartheta}(\theta/\vartheta)^\alpha P_{{\theta}}f_{\Omega}(\alpha
|M)\in E(\widetilde{M}_r)$, where $E(\widetilde{M}_r)$ is the
space of functions $\rho(\alpha) $ of the form
 $$\rho(\alpha)\,\,=
\sum_{\theta\in A(\widetilde{M}_r)\setminus e}
\theta^{\alpha}p_{\theta}(\alpha),\,\,\, p_{\theta}(\alpha)\in
D((\widetilde{M}_r)_{\theta}).
$$
We now consider the sum
\begin{equation}\label{Eq:dtpsumBspline}
\sum_{\alpha \in {\Z}^s}P_{{\vartheta}}f_{\Omega}(\alpha
|M)B(x-\alpha |\widetilde{M}_r), \mbox{ where }x\in
c(\Omega,H_{{\vartheta}}).
\end{equation}
Let $\Psi(y)=P_{{\vartheta}}f_{\Omega}(y |M)B(x-y
|\widetilde{M}_r)$. Then
\begin{eqnarray*}
\widehat{\Psi}(\xi)&=&\int_{\R^s}\Psi(y)\exp(-iy\xi)dy\\
&=&\int_{\R^s}P_{{\vartheta}}f_{\Omega}(y |M)B(x-y
|\widetilde{M}_r)\exp(-iy\xi)dy\\
&=&\int_{\R^s}P_{{\vartheta}}f_{\Omega}(x-t |M)B(t
|\widetilde{M}_r)\exp(-i(x-t)\xi)dt\\
&=&\exp(-ix\xi)\int_{\R^s}P_{{\vartheta}}f_{\Omega}(x-t |M)B(t
|\widetilde{M}_r)\exp(it\xi)dt\\
 &=&\exp({-ix\xi})P_{{\vartheta}}f_{\Omega}(-iD+x
|M)\widehat{B}(-\xi|\widetilde{M}_r).
\end{eqnarray*}
 Taking into account that
(see \cite{DahmenDM})
$$
q(D)\widehat{B}(2\pi \alpha|\widetilde{M}_r)=0, \mbox{ where }
\alpha \in {\Z}^s\setminus 0\,\, \mbox{ and }\,\, q\in
D(\widetilde{M}_r),
$$
 Poisson's
summation formula converts the sum (\ref{Eq:dtpsumBspline}) into
\begin{equation}\label{Eq:dtpsumBspline1}
P_{{\vartheta}}f_{\Omega}(-iD+x|M)\widehat{B}(\cdot|\widetilde{M}_{r})(0),
\end{equation}
since $P_{\vartheta}f_{\Omega}(\cdot +x|M)\in D(M_{{\vartheta}})$
for each fixed $x$ and $ D(M_{{\vartheta}})\subset
D(\widetilde{M}_r)$. Expanding $P_{{\vartheta}}f_{\Omega}$ in a
Taylor series  for each fixed $x$ and noting that
$\widehat{B}(0|\widetilde{M}_{r})=1$, one has
\begin{eqnarray}\label{eq:lema548}
& &P_{{\vartheta}}f_{\Omega}(-iD+x|M)\widehat{B}(\cdot|\widetilde{M}_{r})(0)\nonumber\\
&=&P_{{\vartheta}}f_{\Omega}(x|M)+\!\! \sum_{0<|u|\leq
n-\kappa-s}D^uP_{{\vartheta}}f_{\Omega}(x|M)(-i)^{|u|}D^u\widehat{B}(\cdot|\widetilde{M}_{r})(0)/u!,
\end{eqnarray}
where $x\in c(\Omega,H_{{\vartheta}})$. So, combining
(\ref{eq:le}), (\ref{eq:peq1}), (\ref{eq:le541}),
(\ref{Eq:dtpsumBspline}), (\ref{Eq:dtpsumBspline1}) and
(\ref{eq:lema548}),
 we have
$$
Q_{{\vartheta},r}(x)=P_{{\vartheta}}f_{\Omega}(x|M)+\sum_{|u|=1}^{n-\kappa-s}D^uP_{{\vartheta}}
f_{\Omega}(x|M)(-i)^{|u|}D^u\widehat{B}(\cdot|\widetilde{M}_{r})(0)/u!.
$$
\end{proof}
Finally, we present an explicit formula for $Q_{\vartheta,r}(x)$.
\begin{lemma}\label{le:Qexpli}
When $x\in c(\Omega,H_{{\vartheta}}),$
$$
Q_{{\vartheta},r}(x)=\sum_{0\leq j_1+\cdots+j_\kappa\leq
n-s-\kappa }\, \, \frac{1}{r^\kappa}D_{m_1}^{j_1}\cdots
D_{m_\kappa}^{j_\kappa}T(x|M_{{\vartheta}})\prod_{i=1}^\kappa
\frac{s_{j_i+1}({\vartheta}^{-m_i})}{(j_i+1)!}.
$$
\end{lemma}
\begin{proof}
 When $x\in c(\Omega, H_{{\vartheta} })$, we can see that
$Q_{{\vartheta},r}(x)$ is a polynomial of degree less than
$n-s+1$. Using the Taylor expansion, we obtain
$$
T(x-\sum_{i=1}^\kappa
r_im_i|\widetilde{M}_{r})=T(x|\widetilde{M}_{r})+
\sum_{j=1}^{n-s}\frac{1}{j!}(-\sum_{i=1}^\kappa
r_iD_{m_i})^jT(x|\widetilde{M}_{r}),
$$
where $x\in c(\Omega, H_{{\vartheta} })$. Noting that
$\sum\limits_{r_i=0}\limits^{r-1}{\vartheta}^{-r_im_i}=0$, we have
\begin{equation}\label{eq:qq}
\begin{split}
 & Q_{\vartheta,r}(x)\\
 &=\sum\limits_{0\leq
r_1,\ldots,r_\kappa<r}{\vartheta}^{-\sum_{i=1}^\kappa
r_im_i}T(x-\sum_{i=1}^\kappa m_ir_i|\widetilde{M}_{r})\\
 &= \sum\limits_{0\leq
r_1,\ldots,r_\kappa<r}{\vartheta}^{-\sum_{i=1}^\kappa
r_im_i}\left(T(x|\widetilde{M}_{r})+
\sum_{j=1}^{n-s}\frac{1}{j!}\left(-\sum_{i=1}^\kappa
r_iD_{m_i}\right) ^jT(x|\widetilde{M}_{r})\right)\\
&= \sum\limits_{0\leq
r_1,\ldots,r_\kappa<r}{\vartheta}^{-\sum_{i=1}^\kappa r_im_i}
\sum_{0\leq j_1+\cdots +j_\kappa \leq n-s}
\left(\prod_{i=1}^\kappa
\frac{(-r_iD_{m_i})^{j_i}}{j_i!}\right) T(x|\widetilde{M}_{r})\\
&= \sum_{0\leq j_1+\cdots +j_\kappa \leq n-s}\sum\limits_{0\leq
r_1,\ldots,r_\kappa<r} {\vartheta}^{-\sum_{i=1}^\kappa r_im_i}
\left(\prod_{i=1}^\kappa \frac{(-r_iD_{m_i})^{j_i}}{j_i!}\right)
T(x|\widetilde{M}_{r}).
\end{split}
\end{equation}
 Moreover, if there exists an index $j_h=0$, where $1\leq h\leq \kappa$,
then
\begin{equation}\label{eq:qq1}
 \sum\limits_{0\leq
r_1,\ldots,r_\kappa<r}{\vartheta}^{-\sum_{i=1}^\kappa
r_im_i}\prod_{i=1}^\kappa (-r_iD_{m_i})^{j_i}=\prod_{i=1}^\kappa
\left(\sum_{r_i=0}^{r-1}{\vartheta}^{-r_im_i}(-r_iD_{m_i})^{j_i}\right)\equiv
0,
\end{equation}
since
$\sum_{r_{h}=0}^{r-1}{\vartheta}^{-r_{h}m_{h}}(-r_{h}D_{m_{h}})^{j_{h}}=
\sum_{r_{h}=0}^{r-1}{\vartheta}^{-r_{h}m_{h}} \equiv 0$.

Recall that $s_j(x)=(-1)^j(x+2^jx^2+\cdots+(r-1)^jx^{r-1})$. Based
on (\ref{eq:qq}), (\ref{eq:qq1}) and (\ref{eq:Tre}), we have
\begin{equation*}
\begin{split}
Q_{{\vartheta},r}(x)&=
 \sum\limits_{0\leq
r_1,\ldots,r_\kappa<r}{\vartheta}^{-\sum_{i=1}^\kappa r_im_i}
\sum_{0\leq j_1+\cdots +j_\kappa \leq n-s}\left(\prod_{i=1}^\kappa
\frac{(-r_iD_{m_i})^{j_i}}{j_i!}\right) T(x|\widetilde{M}_{r})\\
&=\sum_{\kappa \leq j_1+\cdots +j_\kappa \leq n-s\atop j_i\geq 1}
\sum\limits_{0\leq r_1,\ldots,r_\kappa<r} \left(\prod_{i=1}^\kappa
{\vartheta}^{- r_im_i}\frac{(-r_iD_{m_i})^{j_i}}{j_i!}\right) T(x|\widetilde{M}_{r})\\
&=\sum_{\kappa \leq j_1+\cdots +j_\kappa \leq n-s\atop j_i\geq
1}\frac{1}{r^\kappa} \sum\limits_{0\leq r_1<r}\cdots
\sum\limits_{0\leq r_\kappa<r} \left(\prod_{i=1}^\kappa
{\vartheta}^{-
r_im_i}\frac{(-r_i)^{j_i}}{j_i!}D_{m_i}^{j_i-1}\right)
T(x|{M}_{\vartheta})\\
&=\sum_{0 \leq j_1+\cdots +j_\kappa \leq n-s-\kappa
}\frac{1}{r^\kappa} \prod_{i=1}^\kappa \left({\sum\limits_{0\leq
r_i<r} \vartheta}^{-
r_im_i}\frac{(-r_i)^{j_i+1}}{(j_i+1)!}D_{m_i}^{j_i}\right)
T(x|{M}_{\vartheta})\\
&= \sum_{0\leq j_1+\cdots+j_\kappa\leq n-s-\kappa
}\frac{1}{r^\kappa}D_{m_1}^{j_1}\cdots
D_{m_\kappa}^{j_\kappa}T(x|M_{{\vartheta}})\prod_{i=1}^\kappa
\frac{s_{j_i+1}({\vartheta}^{-m_i})}{(j_i+1)!},
\end{split}
\end{equation*}
provided that $x\in c(\Omega,H_{{\vartheta}})$.
\end{proof}

We now have all the ingredients for the proof of our main theorem.

\begin{proof}[Proof of Theorem \ref{Th:FormPtheta}]

Since $P_{{\vartheta}}f_{\Omega}$ is a polynomial, it can be
determined by its values on $c(\Omega,H_{{\vartheta}})$ according
to Lemma \ref{Le:asetpro}. We recall (\ref{eq:qform}) that, for
$\mu\in \Z_+$, the homogeneous polynomials
$q_{\mu,r}^{\vartheta}(x)$
  of degree $n-s-\kappa-\mu$  satisfy
$$
q_{ \mu ,r}^{\vartheta}(x)\,\,=\!\! \sum_{ j_1+\cdots+j_\kappa
=\mu }\frac{1}{r^\kappa}D_{m_1}^{j_1}\cdots
D_{m_\kappa}^{j_\kappa}T(x|M_{{\vartheta}})\prod_{i=1}^\kappa
\frac{s_{j_i+1}({\vartheta}^{-m_i})}{(j_i+1)!},\,\,x\in \Omega.
$$
 Based on
Lemma \ref{le:Qexpli}, we have
$$
Q_{\vartheta,r}(x)=\sum_{0\leq \mu\leq
n-s-\kappa}q_{\mu,r}^{\vartheta}(x),\,\,\,  x\in
c(\Omega,H_{{\vartheta}}).
$$
Note that $P_{\vartheta} f_{\Omega}(x|M)$  can be written in  the
form $\sum_{\mu=0}^{n-s-\kappa}p^{\vartheta}_{\mu,\Omega}(x)$. By
Lemma \ref{le:4}, one has
\begin{equation}\label{eq:comp5}
\sum_{\mu=0}^{
n-s-\kappa}q_{\mu,r}^{\vartheta}(x)\,\,=\,\,P_{{\vartheta}}f_{\Omega}(x|M)+\!\!\!
\sum_{|u|=1}^{n-\kappa-s}\!\!\!
D^uP_{{\vartheta}}f_{\Omega}(x|M)(-i)^{|u|}D^u\widehat{B}(0|\widetilde{M}_{r})/u!.
\end{equation}
Comparing the homogeneous polynomials on both sides of
(\ref{eq:comp5}), we arrive at
\begin{equation*}
\begin{split}
p^{\vartheta}_{0,\Omega}(x)&=q_{0,r}^{\vartheta}(x),\\
p^{\vartheta}_{\mu,\Omega}(x)&=q_{\mu,r}^{\vartheta}(x)
-\sum_{j=0}^{\mu-1}\left(\sum_{|u|=\mu-j}D^up^{\vartheta}_{j,\Omega}(x)
(-i)^{|u|}D^u\widehat{B}(0|\widetilde{M}_{r})/u!\right),\,\,
\mu\geq 1.
\end{split}
\end{equation*}
\end{proof}

\section{Discrete  truncated powers associated with  special matrices }
\setcounter{section}{6}

 In this section, we shall present more detailed information
concerning  $t(\cdot |M)$ for some particular matrices $M$. We
first introduce some definitions. Let
$${\mathscr S}_k(M):=\{
Y\subseteq M:\#Y=s+k,\,\, {\rm span}(Y)={\R}^s\}.$$
 In this notation, ${\mathscr B}(M)={\mathscr S}_0(M).$
  If
 ${\rm gcd} \{ |\det(X)| : X\in {\mathscr
B}(Y) \}=1$ for any $Y\in {\mathscr S}_k(M),$  then $M$ is called
a {\it $k$-prime matrix}. In particular, when $M$ is a $1$-prime
matrix, it is also called a {\it pairwise relatively prime
matrix}. When $s=1$, the $k$-prime matrix has the property that
any $k$ integers in  $\{ m_1,\ldots,m_n\}$ have no common factor.

An explicit formula for $t(\cdot |M)$ is presented in
\cite{BeckFrobenius} when $s=1$ and $M$ is a 1-prime matrix.
\begin{theorem}{\rm (\cite{BeckFrobenius})}
Suppose that $M=(a_1,\ldots,a_n)\in {\Z}^n$ with $a_1,\ldots,a_n$
pairwise relatively prime. Then
$$
t(\alpha |M)=R_{-\alpha }(a_1,\ldots,a_n)+(-1)^n\sum_{j=1}^n
\frac{1}{a_j}\sum_{\theta^{a_j}=1\neq \theta}\prod_{k=1,k\neq
j}^n\frac{\theta^{-\alpha }}{\theta^{a_k}-1},
$$
where $R_{-\alpha }(a_1,\ldots,a_n)=-{\rm Res}(F_{-\alpha }(z),z=1)$,
and $F_{-\alpha }(z)=\frac{z^{-\alpha -1}}{(1-z^{a_1})\cdots
(1-z^{a_n})(1-z)}$.
\end{theorem}

Based on Theorem  \ref{Th:FormPtheta}, we can extend the theorem
above to higher dimensions.
\begin{theorem}\label{Th:form}
Under the conditions of Theorem \ref{Theorem:DaMi,truncated power
structure},
$$
f_{\Omega}(\alpha |M)=P_ef_{\Omega}(\alpha |M)+\sum_{{\theta} \in
A(M)\setminus e}{\theta} ^{\alpha
}\frac{1}{|\det(M_{{\theta}})|}\prod_{w\in M\setminus
M_{{\theta}}}\frac{1}{1-{\theta}^{-w}}
1_{{\rm cone}(M_{{\theta}})}(\Omega),
$$
provided that $M$ is a 1-prime matrix, where $P_ef_{\Omega}(\alpha
|M)$ is presented in Theorem \ref{Th:FormP}.
\end{theorem}
\begin{proof}
Select a $\theta \in A(M)\setminus e$. To this end, we need to
prove that $\#M_{\theta}=s$. By (\ref{the term of A(M)}), $\theta$
has the form
\begin{equation*}
\theta ={\exp}(2\pi i \alpha/|\det Y|),
\end{equation*}
where the definitions of $Y=(y_1,\ldots,y_s)$ and $\alpha$ are
identical to those of $(\ref{the term of A(M)})$. Note that
$\theta^y=1$ provided that $y\in Y$, i.e., $Y\subset M_\theta$. We
assert that $M_{\theta}=Y$. Otherwise, if there exists $y_0\in
M\setminus Y$ such that $\theta ^{y_0}=1$, then we can select  an
$s\times s$ square matrix $Y'\subset (y_0,y_1,\ldots,y_s)$ such
that $y_0\in Y'$ and ${\rm span}(Y')={\R}^s$.  Without loss of
generality, we suppose $Y'=(y_1,\ldots,y_{s-1},y_0)$. For each
$y\in Y'$, we have $\theta^y=1$. Hence, there is a
$\beta:=(\beta_1,\ldots,\beta_s)\in \Z^s$ such that $\theta$ is in
the form of
\begin{equation*}
\theta \,\,=\,\,{\exp}(2\pi i \beta/|\det Y'|).
\end{equation*}
 Without loss of generality, we suppose that $\beta_1\neq 0$.
Since $0< \alpha_1\leq |\det Y|-1$ and $0<\beta_1\leq \det(Y')-1$,
we have ${\rm gcd}(\det(Y)$, $\det(Y'))>1$ due to $\alpha_1/|\det
Y|=\beta_1/|\det Y'|$, which contradicts ${\rm
gcd}(\det(Y'),\det(Y))=1$. That is, $M_{\theta}=Y$.

Note that $P_{\theta}f_{\Omega}(\alpha |M)$ is a constant. This
follows from the fact that  the polynomial
$P_{\theta}f_{\Omega}(\alpha |M)\in D(M_{\theta})\subset {\PP}_0.$
By Theorem \ref{Theorem:leadingpart of all },
$$P_{\theta}f_{\Omega}(\alpha |M)=T(\alpha |M_{\theta})\prod_{w\in
M\setminus M_{\theta}}\frac{1}{1-\theta^{-w}},$$
 when $\alpha\in \Omega$. Then, we have
\begin{equation*}
 f_{\Omega}(\alpha |M)=P_ef_{\Omega}(\alpha |M)+\sum_{\theta \in
A(M)\setminus e}\theta ^{\alpha }\prod_{w\in M\setminus
M_{\theta}}\frac{1}{1-\theta^{-w}}T(\alpha |M_{\theta}),\,\,
\alpha \in \Omega\cap \Z^s.
\end{equation*}
It implies that
\begin{equation*}
 f_{\Omega}(\alpha |M)= P_ef_{\Omega}(\alpha |M)+\sum_{\theta \in A(M)\setminus e}\theta
^{\alpha }\frac{1}{|\det(M_{\theta})|}\prod_{w\in M\setminus
M_{\theta}}\frac{1}{1-\theta^{-w}} 1_{{\rm
cone}(M_{\theta})}(\Omega) .
\end{equation*}
\end{proof}

We next show one example, as the application  of Theorem
\ref{Th:form}.

\begin{example}\label{ex:1}
Consider the following partition function
$$p_{\{a,b\}}(n)=\#\{(x,y)\in {\Z}_+^2 : ax+by=n,a,b\in {\Z}_+
\},
$$
where $a$ and $b$ are relatively prime. We shall prove that
$p_{\{a,b\}}(n)=\frac{n}{ab}-\{\frac{b^{-1}n}{a}\}-\{\frac{a^{-1}n}{b}\}+1,$
which is  ``the beautiful formula due to Popoviciu"
\cite{beckbook}, where $\{\frac{t}{a}\}$ is the fractional part of
$\frac{t}{a}$, $b^{-1}$ and $a^{-1}$ denote two integers
satisfying $b^{-1}b\equiv 1\mod  a$, and $a^{-1}a\equiv 1 \mod
b$.

 When $n\geq 0$, $p_{\{ a,b\}}(n)=t(n|(a,b))$. Note that
$T(x|(a,b))=\frac{x_+}{ab}$ and $D^1\widehat{B}(0
|(a,b))=-\frac{i}{2}(a+b).$ By Theorem \ref{Th:form}, when $n\geq
0$, $p_{\{a,b\}}(n)$ equals to
 \begin{eqnarray}\label{eq:pa2}
\frac{n}{ab}+\frac{a+b}{2ab}+\frac{1}{a}\sum_{k=1}^{a-1}\frac{{\exp}({\frac{2\pi
i  nk }{a}})}{1-{\exp}(-{\frac{2\pi i bk }{a}})}+
\frac{1}{b}\sum_{k=1}^{b-1}\frac{{\exp}({\frac{2\pi i nk
}{b}})}{1-{{\exp}(-\frac{2\pi i a k}{b}})}.
\end{eqnarray}
Based on the discrete Fourier analysis, one has (see
\cite{beckbook}, pp.144)
\begin{equation}\label{eq:dis}
-\{\frac{t}{a}\}=\frac{1-a}{2a}+\frac{1}{a}\sum_{k=1}^{a-1}\frac{{\rm
exp}({\frac{2\pi itk}{a} })}{1-{\rm exp}({-\frac{2\pi ik}{a} })},
\end{equation}
where $t,a\in {\Z}.$ According to (\ref{eq:dis}), (\ref{eq:pa2})
can be reduced to
$\frac{n}{ab}-\{\frac{b^{-1}n}{a}\}-\{\frac{a^{-1}n}{b}\}+1$.
\end{example}

\section{ Exact volume of polytopes and the Ehrhart polynomial}

\setcounter{section}{7}
\subsection{Exact volume of polytopes}

A convex polytope $P$ is the convex hull of a finite set of points
in ${\R}^n$. In this section, we shall omit the qualifier
``convex" since we confine our discussion to such polytopes. An
integer polytope is a polytope whose vertices have integer
coordinates. Similarly, a rational polytope is a polytope whose
vertices have rational coordinates. The exact computation of the
volume of $P$ is an important and difficult problem which has
close ties to various mathematical areas
\cite{BeckBirkoffvolume,volume1,Birkoffvolume,volume}. If $P$ is a
$d$-dimensional polytope in ${\R}^n$, then let ${\rm vol}_n(P)$
denote the $d$-dimensional volume of $P$ in ${\R}^n$, i.e., the
$d$-dimensional measure of $P$ in ${\R}^n$. Let ${\R}_P$ denote
the affine space that is spanned by the vertex vectors of $P$. The
lattice points in ${\R}_P$ form an Abelian group of rank $d$,
i.e., ${\R}_P\cap {\Z}^d$ is isomorphic to ${\Z}^d.$ Hence, there
exists an invertible affine linear transformation $T:{\R}_P
\rightarrow {\R}^d$ satisfying $T({\R}_P\cap {\Z}^n)={\Z}^d.$ The
{\it relative volume} of $P,$ denoted as ${\rm vol}(P),$ is just
the $d$-dimensional volume of the image $T(P)\subset {\R}^d$. For
more detailed information about the relative volume, the reader is
referred  to \cite{stanley}.

We next introduce  a   method for computing the relative volume of
polytopes, which depends on the counting function for the integer
points in a polytope. We consider the function of an
integer-valued variable $g$ that describes the number of lattice
points lying inside the dilated polytope $gP$:
$$
L_P(g):=\# \left(gP \bigcap {\Z}^n\right).
$$
In \cite{ehrhart}, Ehrhart inaugurated the study of the general
properties of $L_P(g)$. He proved that for an  $n$-dimensional
polytope $P$, $L_P(g)$ is a polynomial in the positive integer
variable $g$ and that in fact
\begin{eqnarray}\label{Ehrhartp}
L_P(g)\,\,=\,\,{\rm vol}(P)g^n+\frac{1}{2}{\rm vol}(\partial
P)g^{n-1}+\cdots+\chi(P).
\end{eqnarray}
Here, $\chi (P)$ is the Euler characteristic of $P$ and ${\rm
vol}(\partial P)$ is the surface area of $P$ normalized with
respect to the sublattice on each face of $P$. Moreover, the
leading coefficient of $L_P(g)$ is the relative volume of $P$.
Hence, if we obtain the leading coefficient of $L_P(g)$, we can
know the relative volume of $P$. In \cite{BeckBirkoffvolume} and
\cite{Birkoffvolume}, the leading coefficient of $L_P(g)$ was
computed by the interpolation and the residue theorem
respectively.  We can even present a formula for the leading
coefficient of $L_P(g)$ using discrete truncated powers.  Our
result is
\begin{theorem}\label{volumeT}
Suppose that $P_{M,{\bf b}}:=\{ {\bf x}: M{\bf x}= {\bf b},{\bf
x}\in {\R}_+^n\}$ is an $(n-s)$-dimensional rational polytope. Here,
$M$ is an $s\times n$ integer matrix and ${\bf b}$ is an $s$-vector.
Then ${\rm vol}(P_{M,{\bf b }})=C_0\cdot {T({\bf b}|M)}$, where
$C_0={\gcd}\{ |\det(Y)| :Y\in {\mathscr B}(M)\}$.
\end{theorem}
\begin{proof}
To prove the theorem, we only need to prove that
 $${\rm vol}(P_{M,{\bf b }})={\# \{[[M^T))\cap {\Z}^n\}}\cdot {T({\bf
b}|M)}$$
 since $C_0= {\# \{[[M^T))\cap {\Z}^n\}}$ (see
\cite{St3}).
 Since the dimension of $P_{M,{\bf b}}$ is $n-s$, we have ${\bf b}\in
{\rm cone}^\circ(M)$. By (\ref{Eq:truncatedvolume}), ${T({\bf
b}|M)}/{{\rm vol}(P_{M,{\bf b}})}$ is a constant which is
independent  of ${\bf b} \in {\rm cone}^\circ(M)$. For the integer
matrix $M$, there exists an integer vector ${\bf b}_0\in {\rm
cone}^{\circ}(M)$, such that $P_{M,{\bf b}_0}$ is an integer
polytope. We now consider  $L_{P_{M,{{\bf b}_0}}}(g)$. Through the
definition of $L_{P_{M,{{\bf b}_0}}}(g)$ and $t(\cdot |M)$, we
have
$$
L_{P_{M,{{\bf b}_0}}}(g)=t(g{\bf b}_0|M),
$$
with $g{\bf b}_0\in {\rm cone}^{\circ}(M)$ when $g\in {\Z}_+$.
Based on Theorem \ref{Theorem:the leading terms} and Theorem
\ref{Theorem:leadingpart of all }, the leading term of
$L_{P_{M,{{\bf b}_0}}}(g)$ is
$$T(g{\bf b}_0|M)\left(\sum_{\{ \theta \in A(M):M_{\theta}=M\}
}\theta^{g{\bf b}_0}\right) =g^{n-s}T({\bf b}_0|M)\left(\sum_{\{
\theta \in A(M):M_{\theta}=M\} }\theta^{g{\bf b}_0}\right).$$
Since $L_{P_{M,{{\bf b}_0}}}(g)$ is a polynomial, we obtain that
$\sum_{\{ \theta \in A(M):M_{\theta}=M\} }\theta^{g{\bf b}_0}$ is
a constant, i.e., $\theta^{{\bf b}_0}=1$. Hence, we have
$$
\sum_{\{\theta \in A(M):M_{\theta}=M\} }\theta^{g{\bf
b}_0}=\#\{\theta \in A(M):M_{\theta}=M\}=\# \{[[M^T))\cap
{\Z}^n\}.
$$
So, $$\# \{[[M^T))\cap {\Z}^n\} \cdot {T({\bf b}_0|M)} ={{\rm
vol}(P_{M,{{\bf b}_0}})}.$$
 It follows that
$${T({\bf b}|M)}/ {\rm vol}(P_{M,{{\bf b}}})={T({\bf b}_0|M)}/
{\rm vol}(P_{M,{\bf b}_0})=1/ \# \{ [[M^T))\cap {\Z}^n\},$$ when
${\bf b}\in {\rm cone}^\circ (M)$.
 Hence
${\rm vol}(P_{M,{\bf b }})={\# \{[[M^T))\cap {\Z}^n\}}\cdot
{T({\bf b}|M)}$.
\end{proof}

\begin{theorem}\label{Th:npolytopevolume}
Suppose that $P_{A}^{\bf b}=\{{\bf x}: A{\bf x}\leq {\bf b},{\bf
x}\in {\R}_+^n\}$ is a $n$-dimensional polytope, where $A$ is an
$s\times n$ integer matrix and ${\bf b}$ is an $s$-vector. Then
${\rm vol}_n(P_{A}^{\bf b})=T({\bf b}|M)$, where $M=(A,E_{s\times
s})$.
\end{theorem}
\begin{proof}
 Since $P_{A}^{\bf b}$ is an $n$-dimensional polytope in ${\R}^n$,
  the relative volume of $P_{A}^{\bf b}$ is equal to the
volume of $P_{A}^{\bf b}$. Let $P_{A,{\bf b}}=\{({\bf x,y}):A{\bf
x}+E_{s\times s}{\bf y}= {\bf b},{\bf x}\in {\R}_+^n,{\bf y}\in
{\R}_+^s \}$. Since $E_{s\times s}\subset M$, we have ${\rm
gcd}\{|\det(Y)|:Y\in {\mathscr B}(M)\}=1$. If each component of
${\bf b}$ is rational, then $P_{A,{\bf b}}$ is a rational
polytope. Note that ${\rm vol}_n(P_{A}^{\bf b})={\rm
vol}(P_{A}^{\bf b})={\rm vol}(P_{A,{\bf b}}).$ By Theorem
\ref{volumeT}, ${\rm vol}_n(P_{A}^{\bf b})={\rm vol}(P_{A,{\bf
b}})=T({\bf b}|M)$ when $P_{A,{\bf b}}$ is a rational polytope.
Since both $T({\bf b}|M) $ and ${\rm vol}_n(P_{A,{\bf b}})$ are
continuous at ${\bf b}$ and  the real number can be approximated
by rational numbers, one has ${\rm vol}_n(P_{A}^{\bf b})=T({\bf
b}|M)$  for any ${\bf b}\in {\rm cone}^\circ (M)$.
\end{proof}

How can $T({\cdot}|M)$ be computed? In \cite{numeircalT}, an
efficient method for calculating the multivariate truncated power is
presented.

\begin{theorem}{\rm (\cite{numeircalT})}\label{numT}
Let $M$ be an $s\times n$ matrix with columns $m_1,\ldots,m_n\in
{\Z}^s\setminus  0$ such that the origin is not contained in ${\rm
conv}(M).$ For any $\lambda_1,\ldots,\lambda_n\in {\R}$, and
$x=\sum_{j=1}^n\lambda_j m_j,$
\begin{equation}\label{eq:recT}
T(x|M)\,\,=\,\,\frac{1}{n-s}\sum_{j=1}^n\lambda_jT(x|M \setminus
m_j).
\end{equation}
\end{theorem}

Hence, one can compute $T(x|M)$ according to the recurrence
(\ref{eq:recT}). Combining Theorem \ref{Th:npolytopevolume} and
Theorem \ref{numT}, we can present an iterative method for computing
${\rm vol}_n(P_{A}^{\bf b})$. In fact, the computational complexity
is $O(s^{n-s})$. So, when $d:=n-s$ is fixed, i.e., the dimension of
the corresponding polytope is fixed, the run time of the algorithm
based on recurrence (\ref{eq:recT}) is a polynomial function of $s$.
For a detailed description of the algorithm based on
(\ref{eq:recT}), the reader is referred to the \cite{compsimp1} and
\cite{compsimp2}.

In the following, we shall present a simple proof for the volume
formula for a polytope using this method, which is the central
result in \cite{volume}.

Let
$$
\Pi_n(x):=\{ y\in {\R}^n:y_i\geq 0,\, y_1+\cdots+y_i\leq x_1+\cdots
+x_i, \mbox{ for all }1\leq i\leq n\}
$$
for arbitrary $x:=(x_1,\ldots,x_n)$ with $x_i>0$ for all $i$. Let
$V_n(x):={\rm vol}(\Pi_n(x))$. The function $V_n(x)$ is a
homogeneous polynomial of degree $n$ in the variables
$x_1,\ldots,x_n$. This polynomial arises in a variety of
mathematical fields, such as the calculation of probabilities
derived from empirical distribution functions, the study of
parking functions, and plane partitions. In \cite{volume}, an
explicit formula for $V_n(x)$ is presented using a probabilistic
method. Based on (\ref{eq:recT}), we shall show a simple proof for
the explicit formula for $V_n(x)$.
\begin{theorem}{\rm (\cite{volume})}
For each $n=1,2,\ldots,$
$$
V_n(x)=\sum_{k\in K_n}\frac{x^k}{k!} =\frac{1}{n!}\sum_{k\in
K_n}\binom{n}{k}x^k,
$$
where $K_n:=\{ k\in \Z_+^n :  \sum_{i=1}^jk_i \geq j\,\,\mbox{\rm
for all }\,\, 1\leq j\leq n-1\,\, \mbox{\rm and
}\,\,\sum_{i=1}^nk_i=n \}.$
\end{theorem}
\begin{proof}
 Let ${\bf b}=(b_1,\ldots,b_n)^T,$ where $b_i=\sum_{h=1}^ix_h$,
let $I_i$ denote the $i$th column of the $n\times n$ identity
matrix, and let $\widetilde{I}_i=\sum_{h=i}^{n}I_h$. Let
$M_0=(\widetilde{I}_1,\ldots,\widetilde{I}_n)$ and
$M=(M_0,E_{n\times n }).$ Using these notations, we have
$$
\Pi_n(x)=\{ y\in {\R}_+^n: M_0y\leq {\bf b}\}.
$$
By Theorem \ref{Th:npolytopevolume}, $V_n(x)=T({\bf b}|M)$. To
this end,  we present  an explicit formula for $T({\bf b}|M)$.
Note that ${\bf b}=x_1I_1+(x_1+x_2)I_2+\cdots+(\sum_{h=1}^n
x_h)I_n $. By Theorem \ref{numT},
$$T({\bf b}|M)=\frac{1}{n}(x_1T({\bf b}|M\setminus
I_1)+\cdots +(\sum_{h=1}^nx_h) T({\bf b}|M\setminus I_n)).$$
 Note that
\begin{eqnarray*}
{\bf b}&=&x_1I_1+(x_1+x_2)I_2+\cdots+(\sum_{h=1}^nx_h)I_n\\
&=&x_1I_1+\cdots+(\sum_{h=1}^{i-1}x_h)I_{i-1}+(\sum_{h=1}^{i}x_h)\widetilde{I}_{i}+
x_{i+1}I_{i+1}+\cdots+(\sum_{h=i+1}^{n}x_h)I_n\\
&=&x_1I_1+\cdots+(\sum_{h=1}^{i-1}x_h)I_{i-1}+(\sum_{h=1}^{i}x_h)\widetilde{I}_{i}
+x_{i+1}I_{i+1}+\cdots+(\sum_{h=i+1}^{n}x_h)I_n.
\end{eqnarray*}
For any fixed integer $i$, where $1\leq i\leq n$, we have  $T({\bf
b}|M\setminus {I_i\cup \widetilde{I}_i})=0$ since ${\bf b}\notin
{\rm cone}(M\setminus {I_i\cup \widetilde{I}_i})$. Hence,
\begin{eqnarray*}
T({\bf b}|M\setminus I_i)&=&\frac{1}{n-1}(x_1T({\bf b}|M\setminus
I_i\cup I_1)+\cdots+(\sum_{h=1}^{i-1}x_h)T({\bf b}|M\setminus
I_{i}\cup I_{i-1})\\
& &+x_{i+1}({\bf b}|M\setminus I_{i}\cup
I_{i+1})+(\sum_{h=i+1}^nx_{h})T({\bf b}|M\setminus I_{i}\cup
I_n)).
\end{eqnarray*}
Moreover, for any $i>j$, we have
\begin{eqnarray*}
{\bf b}&=&x_1I_1+\cdots
+(\sum_{h=1}^{j-1}x_h)I_{j-1}+(\sum_{h=1}^{j}x_h)\widetilde{I}_j+x_{j+1}I_{j+1}+\cdots+\\
&
&(\sum_{h=j+1}^{i-1}x_h){I}_{i-1}+(\sum_{h=j+1}^{i}x_h)\widetilde{I}_i+x_{i+1}I_{i+1}+\cdots
+(\sum_{h=i+1}^{n}x_h)I_n.
\end{eqnarray*}
Substitute ${\bf b}$ into $T({\bf b}|M\setminus (I_i\cup I_j))$
and then expand $T({\bf b}|M\setminus (I_i\cup I_j))$ by Theorem
\ref{numT}. Continuing the process, we have
$$
V_n(x)\,\,=\,\,T({\bf b }|M)\,\,=\,\,\frac{1}{n!}\sum_{i\in
\sigma_n }x(i),
$$
where $\sigma_n$ is the set of permutations of the set $\{
1,2,\ldots,n\}$, $x(i)=x(i_1,\ldots,i_n)=x(i_1)\cdots x(i_n)$ and
$x(i_k)=x_{j_k+1}+\cdots + x_{i_k}$. Here, $j_k=\max J_k$ (when
$J_k=\emptyset , j_k=0$), and $J_k=\bigcup_{j\leq k-1,i_j\leq i_k}
\{ i_j\}$ for a fixed permutation $i=\left<i_1,\ldots,i_n\right>$.
To this end, we prove
\begin{equation}\label{eq:satanlvo}
\sum_{i\in \sigma_n}\!\!x(i)\,\,=\,\, \sum_{k\in K_n}{n\choose k}
x^k
\end{equation}
 by
induction. When $n=1,$ it is easy to prove that the conclusion
holds. Suppose that the conclusion holds for $n=m$, i.e.,
\begin{eqnarray}\label{hypo}
\sum_{i\in\sigma_m}x(i)\,\,=\,\, \sum_{k\in K_m}{m\choose k} x^k.
\end{eqnarray}
We now consider the case where $n=m+1.$ Note that $x_{m+1}$
appears in $x(i_k)$ if and only if $i_k=m+1.$ Hence, by
(\ref{hypo}), we have
$$
\sum_{i\in \sigma_{m+1}}x(i)\,\,=\,\,\sum_{k\in K_{m+1}}a_kx^k,
$$
where, if $k_{m+1}=0$,
$$a_k= {m\choose
k_1-1,\ldots,k_m}+\cdots+{m\choose k_1,\ldots,k_m-1}={m+1\choose
k},$$
 if
$k_{m+1}=1$, then $a_k=(m+1){m\choose k_1,\ldots,k_m}={m+1\choose
k_1,\ldots,k_m,k_{m+1}}={m+1\choose k}$. Noting $k_{m+1}\leq 1,$
we have
\begin{eqnarray*}
\sum_{i\in\sigma_{m+1}}x(i)\,\,=\,\, \sum_{k\in
K_{m+1}}{m+1\choose k} x^k.
\end{eqnarray*}
Hence, (\ref{eq:satanlvo}) holds when $n=m+1$.
\end{proof}

\subsection{An explicit formula for Ehrhart polynomials}

The explicit formula for $L_P(g)$ is interesting. For any rational
polytope, Ehrhart proved that $L_P(g)$ is a quasipolynomial in
$g.$ Here a quasipolynomial is an expression of the form
$c_n(g)t^g+\cdots +c_0(g),$ where the $c_i(g)$ are periodic
functions in $g$. In (\ref{Ehrhartp}), three coefficients of
$L_P(g)$ are presented. The other coefficients of $L_P(g)$ have
remained a mystery, even for a general lattice 3-simplex, until
rather recently with the work of Pommersheim \cite{pommersheim} in
${\R}^3$, Kantor and Khovanskii \cite{kantor} in ${\R}^4,$ Cappell
and Shaneson \cite{cappell} in ${\R}^n$, and  Diaz and Robins
\cite{Diaz} in ${\R}^n$. In the following theorem, the explicit
formula for $L_P(g)$ is presented in terms of multivariate
truncated powers, where $P$ is a rational polytope.

\begin{theorem}\label{th:lattice}
Let $(n-s)$-dimensional polytopes $P=\{ x\in {\R}_+^n:M x= {\bf b}
\}$, where $M$ is an $s\times n$ integer matrix and ${\bf b}$ is
an integer $s$-vector. We have
$$
L_P(g)=\sum_{j=0}^{n-s}p_{j,\Omega}({\bf b})g^{n-s-j}+\sum_{\theta
\in A(M)\setminus e}\sum_{j=0}^{n-s-(\#M-\#M_{\theta})}(\theta
^{{\bf b} })^gp_{j,\Omega}^\theta({\bf
b})g^{n-s-(\#M-\#M_{\theta})-j},
$$ where $\Omega$ is a fundamental $M$-cone, $p_{j,\Omega}(x)$ and  $p_{j,\Omega}^\theta (x)$
are presented in Theorem \ref{Th:FormP} and Theorem
\ref{Th:FormPtheta} respectively.
\end{theorem}
According to the definition of $L_P(g)$ and $t(\cdot |M),$ we have
$L_{P}(g)=t(g{\bf b}|M).$ By Theorem \ref{Th:FormP}, Theorem
\ref{Th:FormPtheta} and the properties of multivariate truncate
powers, the theorem follows.

\

\bibliographystyle{amsplain}

\end{document}